\documentclass[12pt]{article}
\usepackage{amssymb}
\usepackage{amsmath}
\newcommand{\la}{\lambda}
\newcommand{\lap}{\mbox{$\bigtriangleup$}}

\newcommand{\ra}{{\mbox{$\rightarrow$}}}
\newcommand{\be}{\begin{equation}}
\newcommand{\ee}{\end{equation}}

\newtheorem{mthm}{Theorem}

\newtheorem{lem}{Lemma}[section]

\newtheorem{rem}{Remark}[section]
\begin{document}

\title{Symmetry of Solutions for a Fractional System }

\author{Yan Li\quad Pei Ma \thanks{Corresponding author.}}
\date{\today}
\maketitle
\begin{abstract}
We consider the following  equations:
\begin{equation*}
\left\{\begin{array}{ll}
(-\lap)^{\alpha/2}u(x)=f(v(x)), \\
(-\lap)^{\beta/2}v(x)=g(u(x)),  &x \in R^{n},\\
u,v\geq 0, &x \in R^{n},
\end{array}
\right.
\end{equation*}
for continuous $f, g$  and $\alpha, \beta \in (0,2)$. Under some natural assumptions on $f$ and $g$, by applying the \emph{method of moving
 planes} directly to the system, we obtain symmetry on non-negative solutions without any decay assumption on the solutions
at infinity.

\end{abstract}
\bigskip

{\bf Key words} The fractional Laplacian, narrow region principle, decay at infinity, method of moving planes, radial symmetry.

\bigskip
\section{Introduction}

The famous Lane-Emden equation
\be
\left\{\begin{array}{ll}
-\lap u=u^p, \\
 u>0, & x \in R^n,
 \end{array}
 \right.
\ee
 has been
the central part in the progression of nonlinear analysis
in the past few decades. Such fundamental  results
 as the critical point theory
(also known as the mountain pass
theory), the a priori estimates, the eigenfunction theory and the
Liouville-type theorems have been obtained
(please see \cite{AR}-\cite{PR} and the references therein).

Naturally, one would want to know if these results are valid for the
Lane-Emden system:
\be\label{s19}
\left\{\begin{array}{ll}
-\lap u(x)=v^q(x), \\
-\lap v(x)=u^p(x), &x \in R^n,\\
u,\;v>0,
\end{array}
\right.
\ee
with $p, q>0$. Among which, the problem of
existence and nonexistence
(known as the
Liouville theorem) of solutions has attracted wide attention, but have
not yet been fully answered.
The famous Lane-Emden conjecture states that

\emph{Problem (\ref{s19}) has no classical solutions in the subcritical
case  $\frac{1}{p+1}+\frac{1}{q+1}>\frac{n-2}{n}$.}

The conjecture was validated by Mitidieri \cite{Mi} for radial solutions.
In the critical case  $\frac{1}{p+1}+\frac{1}{q+1}=\frac{n-2}{n}$
and the supercritical case  $\frac{1}{p+1}+\frac{1}{q+1}<\frac{n-2}{n}$,
it has been proved in \cite{SZ} and \cite{Mi} that (\ref{s19}) has
 bounded radial classical solutions.

 In $R^3$, Serrin and Zou \cite{SZ} confirmed the conjecture
 on condition that $(u,v)$ are bounded by polynomials at infinity.
 A decade later,
 Pol$\acute{a}\breve{c}$ik, Quittner and Souplet \cite{PQS} removed the
 growth assumption and obtained same results.
 Recently, still  without growth
 restrictions at infinity, Souplet \cite{S} proved the full conjecture in $R^4$ and
 in a particular category in dimensions higher than 4.
Partial results in higher dimensions were also obtained in certain subregions by Felmer and Figueiredo \cite{FF}, Chen and Li
 \cite{CL2}, Busca and Man$\acute{a}$sevich \cite{BM} and Lin \cite{L}.
The Lane-Emden conjecture was also considered for higher order
elliptic systems:
\be
\left\{\begin{array}{ll}
(-\lap)^m u(x)=v^q(x), \\
(-\lap)^m v(x)=u^p(x), &x \in R^n,\\
u,\;v>0,
\end{array}
\right.
\ee
where $p,q>0$ and $n>1$, $m \in N$. Partial results were obtained
by Mitidieri\cite{Mi}, Yan \cite{Y}, Arthur, Yan and Zhao \cite{AYZ}.

The counterpart of the Lane-Emden equation in the
fractional Laplacian setting has also received a lot of attention.

The fractional
Laplacian $(-\lap)^{\alpha/2}$  is a nonlocal pseudo-differential operator
 defined as
\begin{eqnarray}\label{op}
(-\lap)^{\alpha/2} u(x) &=&
 C_{n,\alpha} P V\int_{\mathbb{R}^n}
 \frac{u(x)-u(z)}{|x-z|^{n+\alpha}} dz\\\nonumber
&=& C_{n,\alpha} \, \lim_{\epsilon \ra 0} \int_{\mathbb{R}^n\setminus
B_{\epsilon}(x)} \frac{u(x)-u(z)}{|x-z|^{n+\alpha}} dz,
\end{eqnarray}
for $\alpha \in (0, 2)$. Here PV stands for the Cauchy principle value.

Let
$$L_{\alpha}=\{u: \mathbb{R}^n\rightarrow \mathbb{R} \mid
\int_{\mathbb{R}^n}\frac{|u(x)|}{1+|x|^{n+\alpha}} \, d x <\infty\}.$$
Then it's easy to see that for $u \in L_{\alpha} \cap C_{loc}^{1,1}(R^n) $,
the integral on the right hand side of
 (\ref{op}) is well defined.
 In this paper, we consider the fractional Laplacian in this space.

The fractional Laplacian has been used to
describe problems emerging from multiple fields. For example,
in the  L\'evy process, it
has been studied as  the
the infinitesimal generator from a perspective of the probability
(see \cite{B}, \cite{V}).
Its nature of non-locality also
fits into
 various research subjects as phase transitions, optimization, flame propagation, finance,
 and so on. Readers who are interested please see
\cite{CRS}, \cite{CT}, \cite{DL}, \cite{ABS}, \cite{G}
and the references therein.

In \cite{CFY}, the authors considered
\begin{equation}\label{s43}
\left\{\begin{array}{ll}
(-\lap)^{\alpha/2}u=u^p,  &x \in R^n_+,\\
u\equiv0, &x \not\in R^{n}_+.
\end{array}
\right.
\end{equation}
They first proved that (\ref{s43}) is equivalent to an integral equation
\begin{equation*}
 u(x)=\int_{R^n_+}G(x,y)u^p(y)dy,
\end{equation*}
where $G(x,y)$ is the Green's function for $(-\lap)^{\alpha/2}$ in
$R^n_+$.
Then by applying the \emph{method of moving planes in integral forms},
they proved the non-existence of positive solutions for
$p\leq\frac{n+\alpha}{n-\alpha}$ without assumptions on the growth
of the solutions.
Using the same method, the authors in \cite{ZCCY}
obtained radial symmetry of positive solutions for a system:
\begin{equation*}
\left\{\begin{array}{ll}
(-\lap)^{\alpha/2}u_i=f_i(u_1,u_2,\cdot\cdot\cdot,u_m), & i=1,2,\cdot\cdot\cdot,m, \\
u_i>0, &x \in R^n,
\end{array}
\right.
\end{equation*}
by investigating its corresponding integral system:
\begin{equation*}
\left\{\begin{array}{ll}
u_i=\int_{R^n}\frac{C_{n,\alpha}}{|x-y|^{n-\alpha}}f_i(u(y)),
& i=1,2,\cdot\cdot\cdot,m, \\
u_i>0, &x \in R^n.
\end{array}
\right.
\end{equation*}
For more details on the \emph{method of moving planes in integral forms},
please see \cite{CL}-\cite{TF}.
In \cite{LM}, the authors
 derived existence and uniqueness
of positive viscosity solutions on a smooth bounded domain
$\Omega\subseteq R^n$:
\begin{equation}
\left\{\begin{array}{ll}
(-\lap)^{\alpha/2}u(x)=v^q(x), \\
(-\lap)^{\alpha/2}v(x)=u^p(x), &x \in \Omega,\\
u=v=0, &x \not\in \Omega,
\end{array}
\right.
\end{equation}
when $pq\neq1$, $p,q>0$ in the supercritical case
 $\frac{1}{p+1}+\frac{1}{q+1}>\frac{n+\alpha}{n}$.
Quite recently, Chen, Li and Li \cite{CLL} developed
 a direct
\emph{method of moving planes} for the fractional
Laplacian. With this, they derived symmetry and nonexistence for:
\begin{equation*}
\left\{\begin{array}{ll}
(-\lap)^{\alpha/2}u(x)=u^p(x),\\
u>0,
\end{array}
\right.
\end{equation*}
in $R^n$ and $R^n_+$.
In \cite{QX}, Quaas and Xia proved a non-existence result for
positive viscosity solutions for:
\begin{equation*}
\left\{\begin{array}{ll}
(-\lap)^{\alpha/2}u(x)=v^q(x), \\
(-\lap)^{\alpha/2}v(x)=u^p(x), &x \in R^{n}_+,\\
u=v=0, &x \not\in R^{n}_+,
\end{array}
\right.
\end{equation*}
with $1<p,q<\frac{n+2\alpha}{n-2\alpha}$.

So far, quite some results have been accumulated about
the fractional system involving the same order of operators. For example,
in
\cite{Y}, Yu considered systems
 of the same fractional order with quite general nonlinearities.
 Using the method of moving planes in integral forms, among which,
 the author
 obtained symmetry of positive solutions. Such results on the system
 have also been proved in \cite{M} and \cite{CHL}. Till now,
 few have been presented dealing with equations involving
different  orders.
In this paper, we introduce a new idea -- the iteration method --
to deal with such problems.

We consider the following  system:
\begin{equation}\label{s1}
\left\{\begin{array}{ll}
(-\lap)^{\alpha/2}u(x)=f(v(x)), \\
(-\lap)^{\beta/2}v(x)=g(u(x)),&x \in R^{n}, \\
u,v\geq 0,
\end{array}
\right.
\end{equation}
with $\alpha, \beta>0$.

First, we use the iteration method to establish the
\emph{maximum principles} for the system.

Let $T_\la$ be a hyperplane in $\mathbb{R}^{n}$. Without loss of generality,
we may assume that
$$T_\la=\{x=(x_1,x') \in \mathbb{R}^{n}\mid x_1=\lambda, \lambda\in \mathbb{R}\}.$$
Let $$x^\la=(2\lambda-x_1, x_2, ..., x_n), \quad \Sigma_\la =\{x \in \mathbb{R}^{n} \mid x_1<\lambda\}.$$

\begin{mthm}[Decay at Infinity]\label{s7}

For $0<\alpha, \beta<2$, assume that
$U\in L_{\alpha} (R^n)\cap C_{loc}^{1,1}(\Omega)$,
 $V\in L_{\beta} (R^n)\cap C_{loc}^{1,1}(\Omega)$, and
 $U$, $V$ are lower semi-continuous on $\bar{\Omega}$. If
\begin{equation}\label{s57}
\left\{\begin{array}{ll}
(-\lap)^{\alpha/2}U(x) +c_1(x)V(x)\geq0, \\
 (-\lap)^{\beta/2}V(x) +c_2(x)U(x)\geq0, &x \in \Omega,\\
U(x), V(x)\geq0,&x \in \Sigma_\la\backslash\Omega,\\
U(x^\la)=-U(x),\\
 V(x^\la)=-V(x), &x \in \Sigma_\la,
\end{array}
\right.
\end{equation}
with
\be
c_1(x)\sim o(\frac{1}{|x|^{\alpha}}),\quad
c_2(x)\sim o(\frac{1}{|x|^{\beta}}),  \quad \mbox{for $|x|$ large,}
\label{s38}
\ee
and
$$c_i(x)<0, \quad i=1,2,$$
then there exists a constant $R>0$ ( depending on $c_i(x)$, but is independent of $U, V$ ) such that if
$$
U(\tilde{x})=\underset{\Omega}{\min}\;U(x)<0,
\quad
V(\bar{x})=\underset{\Omega}{\min}\;V(x)<0,
$$
 then
\be
|\tilde{x}|, |\bar{x}|\leq R.\label{s37}
\ee
\end{mthm}

\begin{mthm}[Narrow Region Principle]\label{s48}

 Let $\Omega$ be a bounded narrow region in $\Sigma_\la$, such that it is contained in $\{x \mid \lambda-l<x_1<\lambda \, \}$ with small $l$. For $0<\alpha, \beta<2$, assume that
$U\in L_{\alpha} (R^n)\cap C_{loc}^{1,1}(\Omega)$,
 $V\in L_{\beta} (R^n)\cap C_{loc}^{1,1}(\Omega)$, and
 $U$, $V$ are lower semi-continuous on $\bar{\Omega}$. If
 $c_i(x)<0$, $i=1,2,$ are bounded from below in $\Omega$  and
\begin{equation}\label{s52}
\left\{\begin{array}{ll}
(-\lap)^{\alpha/2}U(x) +c_1(x)V(x)\geq0, \\
 (-\lap)^{\beta/2}V(x) +c_2(x)U(x)\geq0, &\mbox{ in } \Omega,\\
U(x), V(x)\geq0, &\mbox{ in } \Sigma_\la\backslash\Omega,\\
U(x^\la)=-U(x),\\
 V(x^\la)=-V(x), &\mbox{ in } \Sigma_\la,
\end{array}
\right.
\end{equation}
then for sufficiently small $l$, we have
\begin{equation}\label{s51}
U(x), V(x)  \geq0 \mbox{ in } \Omega.
\end{equation}

If $\Omega$ is unbounded, the conclusion still holds under the condition that
$$\underset{|x| \ra \infty}{\underline{\lim}}U(x), V(x)  \geq0 .$$

Further, if $U(x)$ or $V(x)$ attains 0 somewhere in $\Sigma_\la$, then
\be\label{s80}
U(x)=V(x)\equiv0, \quad x \in R^n.
\ee
\end{mthm}

With the above theorem, we use a \emph{direct method of moving
 planes for the fractional Laplacians}
 \cite{CLL} to show that

\begin{mthm}\label{s45}
Assume that for $r\geq 0$, $f,g $ are nonnegative continuous
functions satisfying:

\emph{(a)} $f(r)$ and $g(r)$ are non-decreasing about $r$;

\emph{(b)} $\frac{f(r)}{r^p}$, $\frac{g(r)}{r^q}$ are bounded
near $r=0$
and
non-increasing with
$p=\frac{n+\alpha}{n-\beta}$ and $q=\frac{n+\beta}{n-\alpha}$.

If $u$ and $v$ are nonnegative solutions for
(\ref{s1}), then
\begin{itemize}
  \item either $u$ and $v$ are constant,
  \item or $f(v)=C_1v^{\frac{n+\alpha}{n-\beta}}$ and
  $g(u)=C_2u^{\frac{n+\beta}{n-\alpha}}$.
\end{itemize}
\end{mthm}

In particular, from Theorem \ref{s45}, we have
\begin{mthm}\label{s66}
Assume $f$ and $g$ satisfy the conditions in Theorem \ref{s45}.
If $u$ and $v$ are nonnegative solutions for
(\ref{s1}), then when $\alpha=\beta$,
\begin{itemize}
  \item either $u$ and $v$ are constant,
  \item or $u(x)=C_1\bigg( \frac{c}{c^2+|x-x_0|^2}\bigg)^{\frac{n-\alpha}{2}}$
      and
        $v(x)=C_2\bigg( \frac{c}{c^2+|x-x_0|^2}\bigg)^{\frac{n-\alpha}{2}}$.
\end{itemize}
\end{mthm}

\begin{rem}
The form of the radial solutions were first obtained in
\cite{CLO2} where the authors classified the positive solutions
of
$$
\left\{\begin{array}{ll}
 u(x)=\int_{R^n}
 \frac{v^{\frac{n+\alpha}{n-\alpha}}(y)}{|x-y|^{n-\alpha}}dy\\
 v(x)=\int_{R^n} \frac{u^{\frac{n+\alpha}{n-\alpha}}(y)}{|x-y|^{n-\alpha}}dy.
 \end{array}
 \right.
 $$
More than a decade later,  the author
in \cite{Y} used results in \cite{CLO2} and
obtained the same explicit expressions for
$$
\left\{\begin{array}{ll}
 u(x)=\int_{R^n} \frac{f(v)(y)}{|x-y|^{n-\alpha}}dy\\
 v(x)=\int_{R^n} \frac{g(u)(y)}{|x-y|^{n-\alpha}}dy.
 \end{array}
 \right.
 $$
\end{rem}

The paper is organized as follows. In Section 2, we verify
Theorem \ref{s7} and \ref{s48}.
To better illustrate the idea, we  first prove Theorem \ref{s45}
under
assumptions on
 the decay rate of solutions at infinity in Section 3.
Then
we complete the proof of  Theorem \ref{s45} in Section 4.
In Section 5, we briefly prove Theorem \ref{s66}. Throughout the
paper, we denote $C$, $C_i$, $i \in N$ as positive constants whose values may
vary from line to line.

\section{Proof of the Maximum Principles}

Our proof of the \emph{maximum principles} for the system
 is inspired by that in \cite{CLL}.
 Let
$\tilde{\Sigma}_\la =\{x \mid x^\la \in \Sigma_\la \}$.

\emph{Proof of Theorem \ref{s7}.}
By (\ref{op}), we have
\begin{eqnarray*}
(-\lap)^{\alpha/2}U(\tilde{x})&=&C_{n,\alpha}PV\int_{\mathbb{R}^n}\frac{U(\tilde{x})-U(y)}
{|\tilde{x}-y|^{n+\alpha}}dy\\
&=&C_{n,\alpha}PV\left\{\int_{\Sigma_\la} \frac{U(\tilde{x})-U(y)}
{|\tilde{x}-y|^{n+\alpha}}dy+\int_{R^n \backslash \Sigma_\la} \frac{U(\tilde{x})-U(y)}
{|\tilde{x}-y|^{n+\alpha}}dy \right\}\\
&=&C_{n,\alpha}PV\left\{\int_{\Sigma_\la} \frac{U(\tilde{x})-U(y)}
{|\tilde{x}-y|^{n+\alpha}}dy+\int_{\Sigma_\la} \frac{U(\tilde{x})-U(y^\la)}
{|\tilde{x}-y^\la|^{n+\alpha}}dy\right\}\\
&=&C_{n,\alpha}PV\left\{\int_{\Sigma_\la} \frac{U(\tilde{x})-U(y)}
{|\tilde{x}-y|^{n+\alpha}}dy+\int_{\Sigma_\la} \frac{U(\tilde{x})+U(y)}
{|\tilde{x}-y^\la|^{n+\alpha}}dy\right\}\\
&\leq&
C_{n,\alpha}\int_{\Sigma_\la}
\left\{\frac{U(\tilde{x})-U(y)}{|\tilde{x}-y^\la|^{n+\alpha}}
+\frac{U(\tilde{x})+U(y)}{|\tilde{x}-y^\la|^{n+\alpha}}\right\}dy\\
&=&C_{n,\alpha}\int_{\Sigma_\la} \frac{2U(\tilde{x})}
{|\tilde{x}-y^\la|^{n+\alpha}}dy.
\end{eqnarray*}
For each fixed $\lambda$, for $\tilde{x} \in \Sigma_\la $ and
$|\tilde{x}|$ sufficiently large, there exists a $C>0$ such that
\begin{eqnarray}\label{s13}
\int_{\Sigma_\la} \frac{1}{|\tilde{x}-y^\la|^{n+\alpha}}dy
\geq\int_{B_{3|\tilde{x}|}(\tilde{x})\backslash
B_{2|\tilde{x}|} (\tilde{x})} \ \frac{1}{|\tilde{x}-y|^{n+\alpha}}dy\sim\frac{C}{|\tilde{x}|^{\alpha}}.
\end{eqnarray}
Hence
\begin{equation}\label{s18}
(-\lap)^{\alpha/2} U(\tilde{x})\leq \frac{CU(\tilde{x})}{|\tilde{x}|^\alpha}<0.
\end{equation}
Together with (\ref{s57}), it's easy to deduce that
\begin{equation}\label{s54}
 V(\tilde{x})<0,
\end{equation}
and
\begin{equation}\label{s14}
  U(\tilde{x})\geq -C
  c_1(\tilde{x})|\tilde{x}|^\alpha V(\tilde{x}).
\end{equation}
From (\ref{s54}), we know that
there exists a $\bar{x}$ such that
$$
V(\bar{x})=\underset{\Omega}{\min}\;V(x)<0.
$$
Similar to  (\ref{s14}), we can derive that

\begin{equation*}
(-\lap)^{\beta/2} V(\bar{x})
 \leq \frac{CV(\bar{x})}{|\bar{x}|^\beta}<0.
\end{equation*}
Combining (\ref{s57}) and  (\ref{s14}), for
$\la$ sufficiently negative, we have
\begin{eqnarray*}
0&\leq&(-\lap)^{\beta/2} V(\bar{x})+ c_2(\bar{x})U(\bar{x})\\
 &\leq&  \frac{CV(\bar{x})}{|\bar{x}|^\beta}
 + c_2(\bar{x})U(\tilde{x})\\
   &\leq& C\bigg( \frac{V(\bar{x})}{|\bar{x}|^\beta}
 -c_2(\bar{x})
  c_1(\tilde{x})|\tilde{x}|^\alpha V(\tilde{x})\bigg)\\
   &\leq&  C\bigg(\frac{V(\bar{x})}{|\bar{x}|^\beta}
 -c_2(\bar{x})
  c_1(\tilde{x})|\tilde{x}|^\alpha V(\bar{x})\bigg)\\
   &\leq&  \frac{CV(\bar{x})}{|\bar{x}|^\beta}
   (1-
  c_1(\tilde{x})|\tilde{x}|^\alpha
  c_2(\bar{x})|\bar{x}|^\beta)\\
  &<&0.
\end{eqnarray*}
The last inequality follows from assumptions
(\ref{s38}).
This contradiction shows that
(\ref{s37})  must be true.
Through an entirely identical argument one can prove the rest of
(\ref{s37}).

This completes the proof.

\medskip

\emph{Proof of Theorem \ref{s48}.}
If (\ref{s51}) does not hold,
then the lower semi-continuity of $U$ on $\bar{\Omega}$ guarantees
that there exists some
$\tilde{x}\in\bar{\Omega}$ such that
$$U(\tilde{x})=\underset{\bar{\Omega}}{\min} \,U<0.$$
And one can further deduce from condition (\ref{s52}) that $\tilde{x}$ is in the interior of $\Omega$.

By (\ref{op}), we have
\begin{eqnarray*}
(-\lap)^{\alpha/2}U(\tilde{x})&=&C_{n,\alpha}PV
\int_{\mathbb{R}^n}\frac{U(\tilde{x})-U(y)}
{|\tilde{x}-y|^{n+\alpha}}dy\\
&\leq& C_{n,\alpha}\int_{\Sigma_\la} \frac{2U(\tilde{x})}
{|\tilde{x}-y^\la|^{n+\alpha}}dy.
\end{eqnarray*}
Let $D=B_{2l}(\tilde{x})\cap \tilde{\Sigma}_\la $. Then we have
\begin{eqnarray}\nonumber
\int_{\Sigma_\la} \frac{1}{|\tilde{x}-y^\la|^{n+\alpha}}dy
&\geq&\int_D \frac{1}{|\tilde{x}-y|^{n+\alpha}}dy\\\nonumber
&\geq & \frac{1}{10}\int_{B_{2l}(\tilde{x})}\frac{1}{|\tilde{x}-y|^{n+\alpha}}dy\\\label{s39}
&\geq &\frac{C}{l^\alpha}.
\end{eqnarray}
Thus,
\begin{equation}\label{s17}
(-\lap)^{\alpha/2}U(\tilde{x})\leq \frac{CU(\tilde{x})}{l^\alpha}<0.
\end{equation}
Together with (\ref{s52}), we have
\begin{equation}\label{s53}
 U(\tilde{x})\geq - c_1(\tilde{x})l^\alpha V(\tilde{x}).
\end{equation}
From (\ref{s53}), we know that
there exists a $\bar{x}$ such that
$$
V(\bar{x})=\underset{\Omega}{\min}\;V(x)<0.
$$
Similar to  (\ref{s17}), we can derive that
\begin{equation*}
(-\lap)^{\beta/2} V(\bar{x})
 \leq \frac{CV(\bar{x})}{l^\beta}<0.
\end{equation*}
Together (\ref{s53}), for
$l$ sufficiently negative, we have
\begin{eqnarray*}
0&\leq&(-\lap)^{\beta/2} V(\bar{x})+ c_2(\bar{x})U(\bar{x})\\
 &\leq&  \frac{CV(\bar{x})}{l^\beta}
 + c_2(\bar{x})U(\tilde{x})\\
   &\leq& C\bigg( \frac{V(\bar{x})}{l^\beta}
 -c_2(\bar{x})
  c_1(\tilde{x})l^\alpha V(\tilde{x})\bigg)\\
   &\leq&  C\bigg(\frac{V(\bar{x})}{l^\beta}
 -c_2(\bar{x})
  c_1(\tilde{x})l^\alpha V(\bar{x})\bigg)\\
   &\leq&  \frac{CV(\bar{x})}{l^\beta}
   (1-c_1(\tilde{x})
  c_2(\bar{x})l^{\alpha+\beta})\\
  &<&0.
\end{eqnarray*}
This contradiction shows that
(\ref{s51})  must be true.  

To prove (\ref{s80}), we suppose that there exists 
$\eta \in \Omega$ such that 
\begin{equation*}
  V(\eta)=0.
\end{equation*}
Then \begin{eqnarray}
&&(-\lap)^{\alpha/2}V(\eta)\nonumber\\
&=&C_{n,\alpha}PV\int_{\mathbb{R}^n}
\frac{-V(y)}
{|\eta-y|^{n+\alpha}}dy\nonumber\\
&=&C_{n,\alpha}PV\int_{\Sigma_\la} \frac{-V(y)}
{|\eta-y|^{n+\alpha}}dy+\int_{\Sigma_\la} \frac{-V(y^\la)}
{|\eta-y^\la|^{n+\alpha}}dy\nonumber\\\label{s81}
&=&C_{n,\alpha}PV\int_{\Sigma_\la} \bigg(\frac{1}
{|x^0-y^\la|^{n+\alpha}}-\frac{1}
{|x^0-y|^{n+\alpha}} \bigg) \,V(y)\,dy.
\end{eqnarray}
If $V(x)\not\equiv0$, then (\ref{s81}) implies that
$$(-\lap)^{\alpha/2}V(\eta)<0.$$
Together with (\ref{s52}), it shows that
$$U(\eta)<0.$$
This is a contradiction with  (\ref{s51}). Hence $ V(x)$ must be identically 0 in $\Sigma_\la$. Since $$V(x^\la)=-V(x),  x \in \Sigma_\la,$$
it shows that 
$$V(x)\equiv0, \quad x \in R^n.$$
Again from (\ref{s52}), we know that
$$U(x)\leq0, \quad x \in\Sigma_\la.$$
Since we already know that
$$U(x)\geq 0, x \in \Sigma_\la,$$
it must hold that
$$U(x)= 0, x \in \Sigma_\la.$$
Together with $U(x^\la)=U(x)$, we arrive at
$$U(x)\equiv 0, x \in R^n.$$

Similarly, one can show that if $U(x)$
attains 0 at one point in $\Sigma_\la$, then both 
$U(x)$ and $V(x)$ are identically 0 in $R^n$.

This completes the proof.

\section{Radial Symmetry under Decay-rate Assumption}

To carry
out the \emph{method of moving planes},  we need to know the behavior of the
solutions at infinity.  We first prove Theorem \ref{s45} in an easier
case by  assuming that for $|x|$ large, $b(p-1)\geq \beta$ and $a(q-1)>\alpha$ (or $b(p-1)> \beta$ and $a(q-1)\geq\alpha$),
\begin{equation*}
u(x)\sim\frac{1}{|x|^a},\quad v(x)\sim\frac{1}{|x|^b},
\end{equation*}
In the next section,  one can see that even without such growth assumption,  Theorem \ref{s45} still holds.

Choose an arbitrary direction for the $x_1$-axis. Let
\begin{equation*} u_\la(x)=u(x^\la),
\end{equation*}
\begin{equation*}
  \tilde{U}_\la(x)= u_\la(x)- u(x),\quad\tilde{ V}_\la(x)=v_\la(x)-v(x).
\end{equation*}

\smallskip
\emph{Step 1. Start moving the plane $T_\la$ from $-\infty$ to the right along the
$x_1$-axis.}
\smallskip

We will show that for $\la$ sufficiently negative,
\begin{equation*}
  \tilde{U}_\la(x) ,  \tilde{V}_\la(x) \geq0, \quad x \in \Sigma_\la \backslash \{  0^\la\}.
\end{equation*}

For fixed $\la$, by the decay rate, we know that
\begin{equation*}
u_\la(x) \: \ra \:0, \quad |x| \:\ra\:\infty.
\end{equation*}
As $\la\:\ra\: -\infty$, for $x \in \Sigma_\la $,
\begin{equation*}
  u(x)\:\ra\:0.
\end{equation*}
Thus  for $x \in \Sigma_\la $,
\begin{equation}\label{s12}
 \tilde{U}_\la(x) \:\ra\:0, \quad \la\:\ra\: -\infty.
\end{equation}
Similarly, one can show that for $x \in \Sigma_\la $,
\begin{equation*}
 \tilde{V}_\la(x) \:\ra\:0, \quad \la\:\ra\: -\infty.
\end{equation*}
If
\begin{equation*}
  \Sigma_{\tilde{U}_\la}^-=\{ x \in \Sigma_\la  \mid \tilde{U}_\la(x)<0 \}\neq \emptyset,
\end{equation*}
then by (\ref{s12}) and Lemma \ref{s46}, we know that there must exist
some $\tilde{x} \in  \Sigma_\la $ such that
\begin{equation*}
 \tilde{ U}_\la(\tilde{x})=\underset{\Sigma_\la}{\min}\:\tilde{ U}_\la <0.
\end{equation*}
On one hand, by (\ref{s18}),
\begin{equation}
(-\lap)^{\alpha/2} \tilde{U}_\la(\tilde{x})\leq \frac{C\tilde{U}_\la(\tilde{x})}{|\tilde{x}|^\alpha}<0.
\end{equation}
On the other hand, by (\ref{s1}),
\begin{equation*}
(-\lap)^{\alpha/2} \tilde{U}_\la(\tilde{x})=f(v_\la(\tilde{x}))-f(v(\tilde{x})).
\end{equation*}
Therefore, by the monotonicity of $f$, we have
\begin{equation*}
 \tilde{V}_\la(\tilde{x})<0.
\end{equation*}
This implies that there exists some
$\bar{x} \in  \Sigma_\la $ such that
\begin{equation*}
 \tilde{ V}_\la(\bar{x})=\underset{\Sigma_\la}{\min}\:\tilde{ V}_\la<0.
\end{equation*}
By the\emph{ mean value theorem}, we have
\begin{eqnarray}
&&(-\lap)^{\alpha/2} \tilde{U}_\la(\tilde{x})\\\nonumber
 &=&f(v_\la(\tilde{x}))-f(v(\tilde{x}))\\\nonumber
 &=& \frac{f(v_\la(\tilde{x}))}{v^p_\la(\tilde{x})}v^p_\la(\tilde{x})-
 \frac{f(v(\tilde{x}))}{v^p(\tilde{x})}v^p(\tilde{x}) \\\nonumber
   &\geq& \frac{f(v(\tilde{x}))}{v^p(\tilde{x})}
   [v^p_\la(\tilde{x})- v^p(\tilde{x})]\\\nonumber
     &=& \frac{f(v(\tilde{x}))}{v^p(\tilde{x})}p \xi^{p-1}(\tilde{x}) \tilde{V}_\la(\tilde{x})
     \qquad \xi \in [ v_\la(\tilde{x}), v(\tilde{x}) ]\\\label{s15}
   &\geq & \frac{f(v(\tilde{x}))}{v^p(\tilde{x})}p v^{p-1}(\tilde{x}) \tilde{V}_\la(\tilde{x}).
\end{eqnarray}
Through a similar argument, one can show that
\begin{equation*}
\tilde{U}_\la(\bar{x})<0,
\end{equation*}
and
$$(-\lap)^{\beta/2} \tilde{V}_\la(\bar{x})
 \geq \frac{g(u(\bar{x}))}{ u^q(\bar{x})}qu^{q-1}
 (\bar{x})\tilde{U}_\la(\bar{x}),
$$
Let $$c_1(x)=-
\frac{f(v(\tilde{x}))}{v^p(\tilde{x})}p v^{p-1}(\tilde{x}),$$
and
$$c_2(x)=-\frac{g(u(\bar{x}))}{ u^q(\bar{x})}qu^{q-1}
 (\bar{x}).$$

By assumption \emph{(b)}, we know that for
$|\tilde{x}|$, $|\bar{x}|$ large,
\begin{equation*}
 c_1(\tilde{x})\sim\frac{1}{|\tilde{x}|^{b(p-1)}}
 \sim o(\frac{1}{|\tilde{x}|^\alpha}),
\end{equation*}
$$c_2(\bar{x})\sim\frac{1}{|\bar{x}|^{a(q-1)}}
 \sim o(\frac{1}{|\bar{x}|^\beta}).$$
By Lemma \ref{s7} (\emph{decay at infinity}), for $\la$ sufficiently negative (less than the $R$ in Lemma \ref{s7}),
it holds that
\begin{equation}\label{s16}
  \tilde{U}_\la(x),  \tilde{V}_\la(x)\geq0, \quad x \in \Sigma_\la.
\end{equation}
This completes step 1.

\smallskip

\emph{Step 2. Continue to move the plane $T_\la$ until the limiting position}
\begin{equation*}
  \la_o=\sup\{ \la<\infty \mid  U_\rho(x), V_\rho(x)\geq0, \:x \in
    \Sigma_\rho, \forall \rho\leq \la \}.
\end{equation*}

Obviously,
\begin{equation*}
   \la_o<\infty.
\end{equation*}
Otherwise, for any $\la>0$,
\begin{equation*}
 u(0^\la) >u(0)\geq0.
\end{equation*}
Meanwhile,
\begin{equation*}
 u(0^\la) \sim \frac{1}{|0^\la|^{n-\alpha}} \: \ra\:0, \quad \la\: \ra\: \infty.
\end{equation*}
This is a contradiction.

Similarly, one can move the plane $T_\la$ from $+\infty$ to the left and show that
\begin{equation*}
  U_{\la_o}(x), V_{\la_o}(x)\leq0, \:x \in
    \Sigma_{\la_o}.
\end{equation*}
Thus,
\begin{equation*}
  U_{\la_o}(x), V_{\la_o}(x)\equiv0, \:x \in
    \Sigma_{\la_o}.
\end{equation*}
Due to the arbitrariness of the $x_1$ direction, we conclude that
$u$, $v$ are symmetric about some point in $R^n$.

\section{Proof of Theorem \ref{s45}}

  Without any decay assumption on $u$ and $v$, we first
 consider their Kelvin transform.
 For any $x^o \in R^n$, let $\bar{u}$ and $\bar{v}$ be the Kelvin transform of $u$ and $v$
 respectively:
 \begin{equation}\label{s35}
   \bar{u}=\frac{1}{|x-x^o|^{n-\alpha}}u\bigg(\frac{x-x^o}{|x-x^o|^2}+x^o \bigg),
 \end{equation}
  \begin{equation}\label{s36}
   \bar{v}=\frac{1}{|x-x^o|^{n-\beta}}v\bigg(\frac{x-x^o}{|x-x^o|^2}+x^o \bigg).
 \end{equation}

  Without loss of generality, let $x^o=0$, then
   \begin{equation*}
   \bar{u}=\frac{1}{|x|^{n-\alpha}}u\bigg(\frac{x}{|x|^2} \bigg),
 \end{equation*}
  \begin{equation}\label{s2}
   \bar{v}=\frac{1}{|x|^{n-\beta}}v\bigg(\frac{x}{|x|^2} \bigg).
 \end{equation}
By (\ref{s1}) and (\ref{s2}),
\begin{eqnarray}
  (-\lap)^{\alpha/2}\bar{u}(x) &=& \frac{1}{|x|^{n+\alpha}} [ (-\lap)^{\alpha/2}u]\bigg(\frac{x}{|x|^2} \bigg)\\
  &=& \frac{1}{|x|^{n+\alpha}}f(|x|^{n-\beta}\bar{v}(x)).
\end{eqnarray}
Similarly,
\begin{equation*}
 (-\lap)^{\beta/2}\bar{v}(x) =\frac{1}{|x|^{n+\beta}}g(|x|^{n-\alpha}\bar{u}(x)).
\end{equation*}

Let
\begin{equation*}
  U_\la(x)= \bar{u}_\la(x)- \bar{u}(x),\quad V_\la(x)=\bar{v}_\la(x)- \bar{v}(x).
\end{equation*}
Then
\begin{equation}\label{s3}
(-\lap)^{\alpha/2} U_\la(x)=\frac{f(|x^\la|^{n-\beta}\bar{v}_\la(x))}{|x^\la|^{n+\alpha}}
-\frac{f(|x|^{n-\beta}\bar{v}(x))}{|x|^{n+\alpha}},
\end{equation}
\begin{equation*}
(-\lap)^{\beta/2} V_\la(x)=\frac{g(|x^\la|^{n-\alpha}\bar{u}_\la(x))}{|x^\la|^{n+\beta}}
-\frac{g(|x|^{n-\alpha}\bar{u}(x))}{|x|^{n+\beta}}.
\end{equation*}

\smallskip
\emph{Step 1. Start moving the plane $T_\la$ from $-\infty$ to the right along the $x_1$-axis.}
\smallskip

We will show that for $\la$ sufficiently negative,
\begin{equation*}
 U_\la(x),  V_\la(x)\geq 0, x \in \Sigma_\la \backslash \{0^\la\}.
\end{equation*}

We noticed that for $\la$ sufficiently negative, there exists a positive constant $C$ such that
\begin{equation*}
 U_\la(x),  V_\la(x)\geq C>0, x \in B_\varepsilon(0^\la) \backslash \{0^\la\}.
\end{equation*}
The proof will be given in the appendix.

Meanwhile, by definition, for $\la$ fixed,
\begin{equation*}
  U_\la(x),  V_\la(x)\: \ra\:0, \mbox{ as } |x| \:\ra \:\infty.
\end{equation*}
Hence if
\begin{equation*}
  \Sigma_{U_\la}^-=\{x \in \Sigma_\la \mid   U_\la(x)<0 \}\neq\emptyset,
\end{equation*}
there must exist a point $\tilde{x}$ such that
\begin{equation*}
  U_\la(\tilde{x})=\underset{\Sigma_\la}{\min}\: U_\la <0.
\end{equation*}
Through arguments similar to those in Lemma \ref{s7}, we obtain
\begin{equation}\label{s4}
(-\lap)^{\alpha/2} U_\la(\tilde{x})\leq \frac{U_\la(\tilde{x})}{|\tilde{x}|^\alpha}<0.
\end{equation}
We claim that
\begin{equation}\label{s5}
  V_\la(\tilde{x})<0.
\end{equation}
Indeed, otherwise, we have $ V_\la(\tilde{x})\geq0$. From (\ref{s3}),
\begin{eqnarray}\nonumber
 &&(-\lap)^{\alpha/2} U_\la(\tilde{x})\\\nonumber
 &=&
 \frac{f(|\tilde{x}^\la|^{n-\beta}\bar{v}_\la(\tilde{x}))}{|\tilde{x}^\la|^{n+\alpha}}
-\frac{f(|\tilde{x}|^{n-\beta}\bar{v}(\tilde{x}))}{|\tilde{x}|^{n+\alpha}} \\\nonumber
   &=&\frac{f(|\tilde{x}^\la|^{n-\beta}\bar{v}_\la(\tilde{x}))}{[|\tilde{x}^\la|^{n-\beta}
   \bar{v}_\la(\tilde{x})]^p }\bar{v}^p_\la(\tilde{x})
-\frac{f(|\tilde{x}|^{n-\beta}\bar{v}_\la(\tilde{x}))}{[|\tilde{x}|^{n-\beta}\bar{v}_\la(\tilde{x})]^p }\bar{v}_\la^p(\tilde{x})\\\label{s9}
   &&+ \frac{f(|\tilde{x}|^{n-\beta}\bar{v}_\la(\tilde{x}))}{|\tilde{x}|^{n+\alpha}} -\frac{f(|\tilde{x}|^{n-\beta}\bar{v}(\tilde{x}))}{|\tilde{x}|^{n+\alpha}}  \\\nonumber
    &\geq&0.
\end{eqnarray}
This is a contradiction with (\ref{s4}), and it proves (\ref{s5}).

From  (\ref{s5}), we can see that there exists a point
$\bar{x}$ such that
\begin{equation*}
  V_\la(\bar{x})=\underset{\Sigma_\la}{\min}\: V_\la(x) <0,
\end{equation*}
and
\begin{eqnarray}\nonumber
  (-\lap)^{\alpha/2} U_\la(\tilde{x}) &=&
  \frac{f(|\tilde{x}^\la|^{n-\beta}\bar{v}_\la(\tilde{x}))}{|\tilde{x}^\la|^{n+\alpha}}
-\frac{f(|\tilde{x}|^{n-\beta}\bar{v}(\tilde{x}))}{|\tilde{x}|^{n+\alpha}}  \\\nonumber
   &=& \frac{f(|\tilde{x}^\la|^{n-\beta}\bar{v}_\la(\tilde{x}))}{[|\tilde{x}^\la|^{n-\beta}
   \bar{v}_\la(\tilde{x})]^p }\bar{v}^p_\la(\tilde{x})
- \frac{f(|\tilde{x}|^{n-\beta}\bar{v}(\tilde{x}))}{[|\tilde{x}^\la|^{n-\beta}
   \bar{v}(\tilde{x})]^p }\bar{v}^p(\tilde{x})\\\nonumber
   &\geq&\frac{f(|\tilde{x}|^{n-\beta}\bar{v}(\tilde{x}))}{[|\tilde{x}^\la|^{n-\beta}
   \bar{v}(\tilde{x})]^p }[\bar{v}^p_\la(\tilde{x})-\bar{v}^p(\tilde{x})]  \\\nonumber
   &=& \frac{f(|\tilde{x}|^{n-\beta}\bar{v}(\tilde{x}))}{[|\tilde{x}^\la|^{n-\beta}
   \bar{v}(\tilde{x})]^p }p\xi^{p-1} V_\la(\tilde{x}),
   \qquad \xi \in [\bar{v}^p_\la(\tilde{x}),\bar{v}^p(\tilde{x})] \\\nonumber
   &\geq& \frac{f(|\tilde{x}|^{n-\beta}\bar{v}(\tilde{x}))}{[|\tilde{x}^\la|^{n-\beta}
   \bar{v}(\tilde{x})]^p }p\bar{v}^{p-1} V_\la(\tilde{x}) \\
   \label{s6}
   &=& \frac{pf(v(\frac{\tilde{x}}{|\tilde{x}|^2}))}
   {|\tilde{x}|^{\alpha+\beta}v(\frac{\tilde{x}}{|\tilde{x}|^2})} V_\la(\tilde{x}).
\end{eqnarray}
Similar to (\ref{s5}), we can show that
$$U_\la(\bar{x})<0.$$
Then we have
\begin{eqnarray*}\nonumber
  (-\lap)^{\beta/2} V_\la(\bar{x}) &=&
  \frac{g(|\bar{x}^\la|^{n-\alpha}\bar{u}_\la(\bar{x}))}{|\bar{x}^\la|^{n+\beta}}
-\frac{g(|\bar{x}|^{n-\alpha}\bar{u}(\bar{x}))}{|\
{x}|^{n+\beta}}  \\\nonumber
   &\geq& \frac{qg(u(\frac{\bar{x}}{|\bar{x}|^2}))}
   {|\bar{x}|^{\alpha+\beta} u(\frac{\bar{x}}{|\bar{x}|^2})} U_\la(\bar{x}).
\end{eqnarray*}
Let
$$c_1(x)=-\frac{pf(v(\frac{\tilde{x}}{|\tilde{x}|^2}))}
{|\tilde{x}|^{\alpha+\beta}v(\frac{\tilde{x}}{|\tilde{x}|^2})},$$
and
$$c_2(x)=-\frac{qg(u(\frac{\bar{x}}{|\bar{x}|^2}))}
   {|\bar{x}|^{\alpha+\beta}u(\frac{\bar{x}}{|\bar{x}|^2})}.$$
From Lemma \ref{s7} (\emph{decay at infinity}), it's easy to deduce
that for $\la$ sufficiently negative,
\begin{equation*}
     U_\la(x),  V_\la(x)\geq0, \quad x \in \Sigma_\la \backslash \{  0^\la\}.
\end{equation*}

\emph{Step 2. Keeping moving the plane $T_\la$ until the limiting position }
\begin{equation*}
  \la_o=\sup\{\la\leq0 \mid  U_\rho(x), V_\rho(x)\geq0, \:x \in
  \Sigma_\rho \backslash \{  0^\rho\}, \forall \rho\leq \la\}.
\end{equation*}

By definition,
\begin{equation}\label{s34}
   U_{\la_o}(x),\: V_{\la_o}(x)\geq0, \quad x \in \Sigma_{\la_o} \backslash \{  0^{\la_o}\}.
\end{equation}

There are two possible cases:
 \begin{itemize}
   \item[i.] either
\begin{equation*}
 U_{\la_o}(x)=V_{\la_o}(x)\equiv0,\quad
 x \in \Sigma_{\la_o}\backslash \{  0^{\la_o}\},
\end{equation*}
   \item[ii.] or
\begin{equation*}
 U_{\la_o}(x), V_{\la_o}(x)>0,\quad
 x \in \Sigma_{\la_o}\backslash \{  0^{\la_o}\}.
\end{equation*}
 \end{itemize}

\emph{Case i. }
Suppose there exists some
$\tilde{x}\in \Sigma_{\la_o}$ such that
$$U_{\la_o}(\tilde{x})=\min_{\Sigma_{\la_o}}U_{\la_o}=0,$$
then it must be true that
\begin{equation}\label{s10}
 U_{\la_o}(x)\equiv0,\quad x \in \Sigma_{\la_o}.
\end{equation}
Otherwise,
\begin{eqnarray}\nonumber
 (-\lap)^{\alpha/2} U_{\la_o}(\tilde{x})
  &=& C_{n,\alpha}PV\int_{\mathbb{R}^n}\frac{-U_{\la_o}(y)}
{|x^0-y|^{n+\alpha}}dy \\\label{s25}
  &< &0.
\end{eqnarray}
On the other hand,
\begin{eqnarray*}\nonumber
 &&(-\lap)^{\alpha/2} U_\la(\tilde{x})\\\nonumber
 &=&
 \frac{f(|\tilde{x}^\la|^{n-\beta}\bar{v}_\la(\tilde{x}))}{|\tilde{x}^\la|^{n+\alpha}}
-\frac{f(|\tilde{x}|^{n-\beta}\bar{v}(\tilde{x}))}{|\tilde{x}|^{n+\alpha}} \\\nonumber
   &=&\frac{f(|\tilde{x}^\la|^{n-\beta}\bar{v}_\la(\tilde{x}))}{[|\tilde{x}^\la|^{n-\beta}
   \bar{v}_\la(\tilde{x})]^p }\bar{v}^p_\la(\tilde{x})
-\frac{f(|\tilde{x}|^{n-\beta}\bar{v}_\la(\tilde{x}))}{[|\tilde{x}|^{n-\beta}\bar{v}_\la(\tilde{x})]^p }\bar{v}_\la^p(\tilde{x})\\
   &&+ \frac{f(|\tilde{x}|^{n-\beta}\bar{v}_\la(\tilde{x}))}{|\tilde{x}|^{n+\alpha}} -\frac{f(|\tilde{x}|^{n-\beta}\bar{v}(\tilde{x}))}{|\tilde{x}|^{n+\alpha}}  \\\nonumber
    &\geq&0.
\end{eqnarray*}
A contradiction with (\ref{s25}).
This proves (\ref{s10}).

Since
\begin{equation*}
 U_{\la_o}(x)=-U_{\la_o}(x^{\la_o}),
\end{equation*}
we have
\begin{equation*}
   U_{\la_o}(x)\equiv0,\quad x \in R^n.
\end{equation*}
Hence
\begin{equation*}
  (-\lap)^{\alpha/2} U_{\la_o}(x)=0.
\end{equation*}
Together with (\ref{s9}), one can deduce that
\begin{equation*}
  \bar{v}_{\la_o}(x)\leq \bar{v}(x),\quad x \in \Sigma_\la.
\end{equation*}
By the definition of $\la_o$,
\begin{equation*}
  \bar{v}_{\la_o}(x)\geq \bar{v}(x),\quad x \in \Sigma_{\la_o}.
\end{equation*}
Thus
\begin{equation*}
V_{\la_o}(x)=\bar{v}_{\la_o}(x)- \bar{v}(x)\equiv0,\quad x \in \Sigma_{\la_o},
\end{equation*}
and
\begin{equation}\label{s11}
V_{\la_o}(x)\equiv0,\quad x \in R^n.
\end{equation}
Similarly, if $V_{\la_o}(x)=0$ somewhere, then we can show that
$$V_{\la_o}(x)=U_{\la_o}(x)\equiv0,\quad x \in R^n.$$
Hence for all $x \in R^n$,
\begin{eqnarray*}
  0 &=&(-\lap)^{\alpha/2} U_{\la_o}(x) \\
   &=&\frac{f(|x^{\la_o}|^{n-\beta}\bar{v}_{\la_o}(x))}{|x^{\la_o}|^{n+\alpha}}
-\frac{f(|x|^{n-\beta}\bar{v}(x))}{|x|^{n+\alpha}}\\
  &=&  \frac{f(|x^{\la_o}|^{n-\beta}\bar{v}_{\la_o}(x))}
  {(|x^{\la_o}|^{n-\beta}\bar{v}_{\la_o}(x))^{\frac{n+\alpha}{n-\beta}}}
  \bar{v}^{\frac{n+\alpha}{n-\beta}}_{\la_o}(x)
  -\frac{f(|x|^{n-\beta}\bar{v}(x))}
  {(|x|^{n-\beta}\bar{v}(x))^{\frac{n+\alpha}{n-\beta}}}
  \bar{v}^{\frac{n+\alpha}{n-\beta}}(x)\\
  &=&\bar{v}^{\frac{n+\alpha}{n-\beta}}(x)
  \bigg(\frac{f(|x^{\la_o}|^{n-\beta}\bar{v}_{\la_o}(x))}
  {(|x^{\la_o}|^{n-\beta}\bar{v}_{\la_o}(x))^{\frac{n+\alpha}{n-\beta}}}
  -\frac{f(|x|^{n-\beta}\bar{v}(x))}
  {(|x|^{n-\beta}\bar{v}(x))^{\frac{n+\alpha}{n-\beta}}} \bigg).
\end{eqnarray*}
Therefore,
\begin{equation*}
 \frac{ f(t)}{t^{\frac{n+\alpha}{n-\beta}}}=C.
\end{equation*}
Similarly, we can prove that
\begin{equation*}
 \frac{ g(t)}{t^{\frac{n+\beta}{n-\alpha}}}=C.
\end{equation*}

\smallskip

\emph{Case ii.} If
\begin{equation}\label{s27}
 U_{\la_o}(x), V_{\la_o}(x)>0, \quad x \in \Sigma_{\la_o},
\end{equation}
then we claim that
\begin{equation}\label{s8}
   \la_o=0.
\end{equation}
To prove (\ref{s8}), suppose that $\la_o<0$.
Then with (\ref{s27}), we would be able to keep moving the plane $T_\la$.
To be precise, for some $\varepsilon>0$ small such
$\la_o+\varepsilon<0$,
it holds that
\begin{equation}\label{s28}
 U_\la(x), V_\la(x)\geq0,
 \quad x \in \Sigma_\la \backslash \{  0^\la\},
 \; \forall\la \in (\la_o, \la_o+\varepsilon).
\end{equation}
This is a contradiction with the definition of
$\la_o$. Therefore,
$$\la_o=0.$$

From (\ref{s27}), we have
\begin{equation*}
 U_{\la_o}(x)>0, \quad
 x \in ( \Sigma_{\la_o} \backslash \{  0^{\la_o}\})\cap B_R(0).
\end{equation*}
Later, in  Lemma \ref{s47}, we will show that
\begin{equation*}
 U_{\la_o}(x)\geq C>0, \quad
 x \in  B_\varepsilon(0^{\la_o})\backslash \{ 0^{\la_o} \}.
 \end{equation*}
Together with bounded-away-from-0 results, we derive that for $\delta>0$ small there exists
some $C$ that
\begin{equation*}
U_{\la_o}(x)\geq C>0, \quad
 x \in (\Sigma_{\la_o-\delta}  \backslash \{  0^{\la_o}\})\cap B_R(0).
\end{equation*}
For $\delta,\varepsilon\ll |\la_o|$,
$$0^\la \in (\Sigma_{\la_o-\delta}  \backslash \{  0^{\la_o}\})
\cap B_R(0).$$
Since $ U_\la(x) $ depends on $\la$ continuously,
\begin{equation}\label{s33}
 U_\la(x) \geq0, \quad
 x \in (\Sigma_{\la_o-\delta}  \backslash \{  0^{\la_o}\})\cap B_R(0).
\end{equation}
What remains is to show that
\begin{equation}\label{s29}
 U_\la(x) \geq0, \quad
 x \in (\Sigma_\la \backslash \Sigma_{\la_o-\delta})\cap B_R(0).
\end{equation}
We argue by contradiction. Suppose (\ref{s29}) does not holds,
then there exists some
$\tilde{x} \in (\Sigma_\la \backslash \Sigma_{\la_o-\delta})\cap
B_R(0)$ such that
\begin{equation*}
U_\la(\tilde{x})=\underset{\Sigma_\la}{\min}\:U_\la <0.
\end{equation*}
Similar to (\ref{s5}), we can show that
\begin{equation*}
  V_\la(\tilde{x})<0.
\end{equation*}
Therefore, there exists some
 $\bar{x} \in  (\Sigma_\la \backslash \Sigma_{\la_o}) \cap
B_R(0)$
such that
\begin{equation*}
V_\la(\bar{x})=\underset{\Sigma_\la}{\min}\:V_\la <0.
\end{equation*}

By (\ref{s39}),
\begin{equation*}
(-\lap)^{\alpha/2} U_\la(\tilde{x})\leq \frac{U_\la(\tilde{x})}{(\delta+\varepsilon)^\alpha}<0.
\end{equation*}
Similarly, one can show that
$$(-\lap)^{\beta/2} V_\la(\bar{x})\leq
\frac{V_\la(\bar{x})}{(\delta+\varepsilon)^\beta}<0.$$

Since
$$\tilde{x}, \bar{x} \in (\Sigma_\la \backslash \Sigma_{\la_o-\delta})\cap
B_R(0),$$
we can choose $\delta,\varepsilon\ll\frac{\la_o}{2}$ such that
$|\tilde{x}|, |\bar{x}|>\frac{2\la_o}{3} $.
Then for some $C$, it holds that
$$v(\frac{\tilde{x}}{|\tilde{x}|^2}),\;
   u(\frac{\bar{x}}{|\bar{x}|^2})\leq C,$$
and
\begin{equation*}
  \frac{1}{ (\delta+\varepsilon)^{\alpha+\beta}}
   \:>\: C+2\:> \:c_1(\tilde{x})c_2(\bar{x}).
\end{equation*}
It then follows from Lemma \ref{s48}(\emph{narrow region principle})
that (\ref{s29}) must be true.

Combining (\ref{s33}) and (\ref{s29}),
we arrive at
\begin{equation*}
 U_\la(x)\geq0,
 \quad x \in \Sigma_\la \backslash \{  0^\la\}.
\end{equation*}
The proof for $V_\la(x)$ is almost the same and we omit it.
This completes the proof of (\ref{s28}) and (\ref{s8}).

Through an identical argument, one can move
$T_\la$ from $+\infty$ to the left and show that
\begin{equation*}
   U_{\la_o}(x),\: V_{\la_o}(x)\leq0, \quad x \in \Sigma_{\la_o} \backslash \{  0^{\la_o}\},
\end{equation*}
with
$$\la_o=0.$$
Together with (\ref{s34}), we obtain
\begin{equation*}
   U_0(x),\: V_0(x)\equiv0, \quad x \in \Sigma_0.
\end{equation*}

For a more general Kelvin transform as in
(\ref{s35}) and (\ref{s36}), through a similar  argument as above
one can show that
$$\la_o=x^o_1,$$ and
\begin{equation*}
   U_{\la_o}(x),\: V_{\la_o}(x)\equiv0, \quad x \in \Sigma_{\la_o}.
\end{equation*}

For any $x^1, x^2 \in R^n$, let their midpoint be the center of
the Kelvin transform:
\begin{equation*}
 x^o=\frac{x^1+ x^2}{2}.
\end{equation*}
Let
\begin{equation*}
  y^1=\frac{x^1-x^o}{|x^1-x^o|^2}+x^o, \quad
  y^2=\frac{x^2-x^o}{|x^2-x^o|^2}+x^o.
\end{equation*}
Then
\begin{equation*}
    y^2=(y^1)^{\la_o},
\end{equation*}
and
\begin{equation*}
  \bar{u}(y^1)=\bar{u}(y^2), \quad\bar{v}(y^1)=\bar{v}(y^2).
\end{equation*}
It thus implies that
$$u(x^1)=u(x^2), \quad v(x^1)=v(x^2).$$
Since  $x^1, x^2$ are arbitrary, $u$ and $v$ must be constant.

\smallskip

Combining \emph{Case i} and \emph{Case ii}, we complete the proof.

\section{Proof of Theorem \ref{s66}}

From Theorem \ref{s45}, we know that when $\alpha=\beta$, if $u$ and $v$ are nonnegative solutions for
(\ref{s1}), then
\begin{itemize}
  \item either $u$ and $v$ are constant,
  \item or $f(v)=C_1v^{\frac{n+\alpha}{n-\alpha}}$ and
  $g(u)=C_2u^{\frac{n+\alpha}{n-\alpha}}$.
\end{itemize}
By Theorem 2 in \cite{ZCCY}, we know that when $\alpha=\beta$,
system (\ref{s1}) is equivalent to
 $$
\left\{\begin{array}{ll}
 u(x)=\int_{R^n}
 \frac{C_1v^{\frac{n+\alpha}{n-\alpha}}(y)}{|x-y|^{n-\alpha}}dy\\
 v(x)=\int_{R^n} \frac{C_2u^{\frac{n+\alpha}{n-\alpha}}(y)}{|x-y|^{n-\alpha}}dy.
 \end{array}
 \right.
 $$
From the results in \cite{CLO2}, we have
 $$u(x)=C_1\bigg( \frac{c}{c^2+|x-x_0|^2}\bigg)^{\frac{n-\alpha}{2}}, \quad
        v(x)=C_2\bigg( \frac{c}{c^2+|x-x_0|^2}\bigg)^{\frac{n-\alpha}{2}}.$$

This proves the theorem.

\section{Appendix}

We use the ideas in the proof of the bounded-away-from-0 lemmas
in the appendix in \cite{CLL}.

\begin{lem}\label{s46}
For $\la$ sufficiently negative, there exists a positive constant $C$ such that
\begin{equation}\label{s20}
 U_\la(x),  V_\la(x)\geq C>0, x \in B_\varepsilon(0^\la) \backslash \{0^\la\}.
\end{equation}
\end{lem}

\textbf{Proof.}  For $\la \in \Sigma_\la $, as
$\la  \:\ra \: -\infty$, it's easy to see that
\begin{equation}\label{s22}
 \bar{u}(x)  \:\ra \: 0.
\end{equation}
To prove (\ref{s20}), it suffices to show that
\begin{equation*}
  \bar{u}_\la(x) \geq C>0, \quad x \in B_\varepsilon(0^\la) \backslash \{0^\la\}.
\end{equation*}
Or equivalently,
\begin{equation*}
  \bar{u}(x) \geq C>0, \quad x \in B_\varepsilon(0) \backslash \{0\}.
\end{equation*}

Let $\eta$ be a smooth
cutoff function such that $\eta(x)\in [0,1]$
in $\mathbb{R}^n$, $supp\,\eta \subset B_2$ and $\eta(x)\equiv1$ in
$B_1$. Let
\begin{equation*}
(-\lap)^{\alpha/2}\varphi(x)=\eta(x)f(v)(x).
\end{equation*}
Then
$$\varphi(x)=C_{n, -\alpha}\int_{\mathbb{R}^n}\frac{\eta(y)f(v)(y)}
{|x-y|^{n-\alpha}}dy=C_{n, -\alpha}\int_{B_2(0)}\frac{\eta(y)f(v)(y)}
{|x-y|^{n-\alpha}}dy.$$
It's trivial that for $|x|$ sufficiently large,
\begin{equation}\label{s21}
  \varphi(x)\sim \frac{1}{|x|^{n-\alpha}}.
\end{equation}
Since
\begin{equation}
\left\{\begin{array}{ll}
(-\lap)^{\alpha/2} (u-\varphi)\geq0,  & x \in B_R,\\
(u-\varphi)(x)\geq0, & x \in B_R^c,
\end{array}
\right.
\end{equation}
by the \emph{maximum principle} (see \cite{Si}), we have
\begin{equation*}
(u-\varphi)(x)\geq0,\quad x \in B_R.
\end{equation*}
Thus
\begin{equation*}
(u-\varphi)(x)\geq0, \quad x\in R^n.
\end{equation*}
For $|x|$ sufficiently large, from (\ref{s21}), one can see that
for some constant $C>0$,
\begin{equation}
u(x) \geq \frac{C}{|x|^{n-\alpha}}.
\end{equation}
Hence for $|x|$ small
\begin{equation*}
u(\frac{x}{|x|^2} ) \geq C|x|^{n-\alpha},
\end{equation*}
and
\begin{equation*}
\bar{u}(x)=\frac{1}{|x|^{n-\alpha}}u(\frac{x}{|x|^2} ) \geq C.
\end{equation*}
Together with (\ref{s22}), it yields that
\begin{equation}\label{s23}
U_\la(x)\geq \frac{C}{2} > 0, \quad x \in B_\varepsilon(0^\la)\backslash\{0^\la\}.
\end{equation}
Through an identical argument, one can show that (\ref{s23})
holds for $V_\la(x)$ as well.

This completes the proof of the lemma.

\begin{lem}\label{s47}
For $\la_o<0$, if either of $U_{\la_o}(x), V_{\la_o}(x) $ is not identically 0, then there exists some
constant $C$ and $\varepsilon>0$ small such that
\begin{equation*}
 U_{\la_o}(x), V_{\la_o}(x)\geq C>0, \quad
 x \in B_\varepsilon(0^{\la_o}) \backslash \{  0^{\la_o}\}.
\end{equation*}
\end{lem}

\textbf{Proof.}
From  Lemma 2.2 in \cite{CFY}, we have the integral expression of
$ U_{\la_o}$:

\begin{eqnarray*}
 U_{\la_o}(x) &=& C_{n, \alpha}\int_{\Sigma_{\la_o}} \bigg(\frac{1}{|x-y|^{n-\alpha}}-
  \frac{1}{|x-y^{\la_o}|^{n-\alpha} }\bigg)\\
  &&\cdot\:\bigg(\frac{f(y^{\la_o}|^{n-\beta}\bar{v}_{\la_o}(y))}{|y^{\la_o}|^{n+\alpha}}
  -\frac{f(|y|^{n-\beta}\bar{v}(y))}{|y|^{n+\alpha}}  \bigg)dy \\
 &=& C_{n, \alpha}\int_{\Sigma_{\la_o}} \bigg(\frac{1}{|x-y|^{n-\alpha}}-
  \frac{1}{|x-y^{\la_o}|^{n-\alpha} }\bigg)\\
  &&\cdot  \: \bigg(\frac{f(|y^{\la_o}|^{n-\beta}\bar{v}_{\la_o}(y))}{[|y^{\la_o}|^{n-\beta}\bar{v}_{\la_o}(y)]^p}
  \:\bar{v}^p_{\la_o}(y)
 -\frac{f(|y|^{n-\beta}\bar{v}_{\la_o}(y))}{[|y|^{n-\beta}\bar{v}_{\la_o}(y)]^p}
  \:\bar{v}^p_{\la_o}(y)\\
  &&+\frac{f(|y|^{n-\beta}\bar{v}_{\la_o}(y))-f(|y|^{n-\beta}\bar{v}(y))}{|y|^{n+\alpha}} \bigg)dy\\
 &\geq& C_{n, \alpha}\int_{\Sigma_{\la_o}} \bigg(\frac{1}{|x-y|^{n-\alpha}}-
  \frac{1}{|x-y^{\la_o}|^{n-\alpha} }\bigg)\\
  &&\cdot\:
 \frac{f(|y|^{n-\beta}\bar{v}_{\la_o}(y))-f(|y|^{n-\beta}\bar{v}(y))}{|y|^{n+\alpha}} dy.
\end{eqnarray*}
Since $$V_{\la_0}(x)\not\equiv 0, \quad x \in \Sigma_{\la_0},$$
there exists some $x^0$ such that
$$V_{\la_o}(x^0)> 0.$$
Thus, for some $\delta>0$ small, it holds that
\begin{equation*}
f(|y|^{n-\beta}\bar{v}_{\la_o}(y))-f(|y|^{n-\beta}\bar{v}(y))\geq C> 0, \quad y \in B_\delta(x^0).
\end{equation*}
Therefore,
\begin{equation}\label{s24}
 U_{\la_o}(x) \geq \int_{ B_\delta(x^0)} C \,dy\geq C>0.
\end{equation}
In a same way, one can show that $ V_{\la_o}(x) $ also satisfies (\ref{s24}).
This completes the proof.

\bigskip

{\em Authors' Addresses and E-mails:}
\medskip

Yan Li

Department of Mathematical Sciences

Yeshiva University

New York, NY, 10033 USA

yali3@mail.yu.edu

\medskip

Pei Ma

Jiangsu Key Laboratory for NSLSCS

School of Mathematical Science

Nanjing Normal University

Nanjing, Jiangsu, 210023 China

Department of Mathematical Sciences

Yeshiva University

New York, NY, 10033 USA

mapei0620@126.com

\end{document}